\documentclass[12pt]{article}
\usepackage{amsmath}
\usepackage{amsthm}
\usepackage{amssymb}
\usepackage{latexsym}
\newtheorem{thm}{Theorem}[section]
\newtheorem{lem}[thm]{Lemma}
\newtheorem{prp}[thm]{Proposition}
\newtheorem{cor}[thm]{Corollary}

\newtheorem{rem}[thm]{Remark}

\newcommand{\nc}{\newcommand*}
\newcommand{\rnc}{\renewcommand*}
\nc\noi{\noindent}
\nc\nn{\nonumber}
\nc{\reff}[1]{(\ref{#1})}
\nc{\ts}{\textstyle}
\nc{\ds}{\displaystyle}
\nc{\mto}{\mapsto}
\nc{\lar}{\leftarrow}
\nc{\rar}{\rightarrow}
\nc{\Lar}{\Leftarrow}
\nc{\Rar}{\Rightarrow}
\nc{\lrar}{\longrightarrow}
\nc{\lra}{\leftrightarrow}
\nc{\hra}{\hookrightarrow}
\nc{\und}{\underline}
\nc{\ov}{\overline}
\nc{\ora}{\overrightarrow}
\nc{\dar}{\downarrow}
\nc{\wt}{\widetilde}
\nc{\mcal}{\mathcal}
\nc{\Cal}{\mathcal}
\rnc{\cal}{\mathcal}
\nc{\goth}{\mathfrak}
\rnc{\bold}{\mathbf}
\rnc{\frak}{\mathfrak}
\rnc{\Bbb}{\mathbb}
\nc{\mbf}[1]{\mbox{\boldmath ${#1}$}}
\nc{\eps}{\epsilon}
\nc{\vep}{\varepsilon}
\nc{\vp}{\varphi}
\nc{\vt}{\vartheta}
\nc{\vr}{\varrho}

\nc{\cA}{\mathcal{A}}
\nc{\cF}{\mathcal{F}}
\nc{\cH}{\mathcal{H}}
\nc{\cP}{\mathcal{P}}
\nc{\cR}{\mathcal{R}}
\nc{\cT}{\mathcal{T}}
\nc{\half}{\frac{1}{2}}
\def\ot{\otimes}
\begin {document}
$$ {} $$
\centerline{\large\bf A.A.Stolin, P.P.Kulish, E.V.Damaskinsky}
\bigskip

\centerline{\Large\bf On construction of universal twist element}
\medskip

\centerline{\Large\bf from $R$-matrix}
\bigskip

\centerline{\bf Abstract}
\medskip

A method to construct the universal twist element
using the constant quasiclassical unitary matrix
solution of the Yang - Baxter equation is proposed.
The method is applied to few known $R$-matrices,
corresponding to Lie (super) algebras of rank one.
\medskip
\bigskip

\section {Introduction}
Quantum groups, as an important class of Hopf algebras, were introduced
by  V.G.Drinfel'd~\cite{D1}. In the theory of quantum groups
the universal $R$-matrix is an essential object intertwining
the coproduct $\Delta$ with the opposite coproduct $\Delta^{\text{op}}$
\begin{equation}\label{a1}
\cR\Delta=\Delta^{\text{op}}\cR.
\end{equation}
Another important object, introduced later, is a twist transformation
(or twist)~\cite{D2,D3}
\begin{equation}\label{a2}
\Delta\to\Delta_t=\cF\Delta\cF^{-1} ,
\end{equation}
where $\cR$ and $\cF$ are some elements of the tensor square of
the Hopf algebra. Twist elements are known explicitly for some quantum
groups, in particular, for quantum deformations of the universal
enveloping algebras of some Lie algebras.

The FRT-formalism~\cite{FRT} permits to reconstruct the quantum group
corresponding to some known matrix representation of universal $R$-matrix.
However this method does not give explicit expressions for the universal
elements.

A similar result for universal twists was obtained by Drinfel'd~\cite{D2}.
It was shown that, if a unitary matrix $R$ being a quasiclassical
solution of the Yang - Baxter equation, is known then there exists a
universal twist $\cF,$ such that $\cR=\cF_{21}{\cF}^{-1}$ in the vector
representation. This situation is similar to the mentioned above:
if a matrix solution of the Yang - Baxter equation is known then
a formal reconstruction of the universal element is possible.
However, as in the previous case, there is no computational method to
obtain explicit expressions for the universal twist $\cF$.
Constructing of universal twist elements is a rather complicated problem.
Some explicit expressions  can be found in papers~\cite{KLMjmp,KLS}
(and refs therein).

In the present paper we suggest a method of constructing of a universal
twist element for a given unitary quasiclassical matrix solutions of the
Yang - Baxter equation. This method is based on an investigation of the
associative algebra defined by the matrix relation
$RT_1 T_2 =T_2 T_1,$ with $T=\left\{t_{ij} \in U(gl_n)^*\right\}.$
This algebra is a deformation of the algebra of polynomial functions on
the group $GL_n$ (eventually, in examples, it is possible to consider its
subgroups). The mentioned above associative algebra is
isomorphic, as a left $U(g\ell_n)$-module algebra, to the nondeformed
(commutative) algebra $Fun$ of polynomials on $GL_n$ (certainly this
isomorphism is not an algebra isomorphism). Let us recall that
there exists a non-degenerate pairing between $Fun$ and $U(g\ell_n)$.
Hence, there is a non-degenerate pairing between $U(g\ell_n)$
and the algebra, generated by the relation  $RT_1 T_2 =T_2 T_1$.
This duality generates a coassociative operation
$d: U(g\ell_n) \to U(g\ell_n)\otimes U(g\ell_n)$
(which is not an algebra homomorphism). Then it turns out that
$\mathcal{F}=d(1)$ is the required twist.

This twist $\cF\in\cH\ot\cH$  allows to introduce in the
Hopf algebra $\cH$ an additional coalgebra structure:
$$
d: \cH\to\cH\ot\cH, \quad d\,:=\,\cF\Delta.
$$
Structures of this type were considered in papers~\cite{Parm,Zakr}
devoted to  analysis of Lie - Poisson groups.

In section 2 we discuss a general connection between the universal twist
$\mathcal{F}$ of a Hopf algebra $\cH$ and the matrix form of the
$RTT$ - relation, which defines a new multiplication law on the dual Hopf
algebra $\cH^*$.
In section 3 we describe some examples of application of this connection
for reconstructing of the universal twisting element from a given $R$-matrix
in the case of Lie (super)algebras of  rank 1. Several details of
calculations are given in the Appendix.

\section {Module algebras and twists}
Let us consider in a given Hopf algebra $\cH$ an element $\cF\in\cH\ot\cH,$
satisfying the twist equation  (the 2-cocycle condition)
$$
\cF_{12}(\Delta\ot id)\cF=\cF_{23}(id\ot\Delta)\cF
$$
and relation
$$(\varepsilon\ot id)\cF=(id\ot\varepsilon)\cF=1\ot 1.$$
Such an element $\cF$ is called the {\it universal twist} or the
{\it universal twist element}.

Define a map  $d_\cF:\cH\rightarrow\cH\ot\cH$,  such that
$d_{\cF}(a)=\cF\Delta (a).$ It is easy to show, that the map $d_{\cF}$
satisfy the following conditions
$$
(d\ot id) d = (id\ot d) d\qquad\text{coassociativity}
$$
and
$$
d_{\cF}(ab)=d_{\cF}(a)\Delta (b).
$$

\begin{lem}
For the twist $\cF$ the following conditions hold:
\begin {enumerate}
\item[\rm{(a)}] $d_{\cF}$ is a coassociative map.
\item[\rm{(b)}] $d_{\cF}(ab)=d_{\cF}(a) \Delta (b)$
\item[\rm{(c)}] $(\varepsilon\ot id) d_{\cF}(a)=1\ot a\quad$ and \quad
$(id\ot\varepsilon) d_{\cF}(a)=a\ot 1$.
\end{enumerate}
Conversely, if $d:\cH\to\cH\ot\cH$ satisfies the conditions
{\rm (a), (b)} and {\rm (c)}, then ${\cF}=d(1)$ is a twist.
\end{lem}

\begin{proof}
We have already shown, that the twist $\cF$ generates a map
$d_{\cF}(a)=\cF\Delta(a),$ satisfying the conditions (a), (b),
while (c) is straightforward.

Now we prove that if the map $d:\cH\to\cH\ot\cH$ is given, then
from the properties (a), (b) it follows, that the map
$$
\cF=d(1)\,:=\,\sum f_1 \ot f_2
$$
satisfies the relation
$$
{\cF}_{12} (\Delta\ot id)(\cF)={\cF}_{23}(id\ot\Delta)({\cF})\,.
$$
Then, by putting $a=1$ in (c) we obtain
$(\varepsilon\ot id)\cF=(id\ot\varepsilon)\cF=1\ot 1$.
\end{proof}

Let us consider the vector space $\cH^\ast$ dual to $\cH$.
We note that $\cH^\ast$
is a right $\cH_{op}$-module since
$$
\langle \alpha, ab \rangle = \langle \alpha_{(1)}, a \rangle
\langle \alpha_{(2)}, b \rangle := \langle \alpha \tilde b, a \rangle,
\quad  \alpha \tilde b = \alpha_{(1)} \langle \alpha_{(2)}, b \rangle.
$$
where $\alpha \in {\cH}^\ast,$  $a,b\in\cH$ and $\cH_{op}$ is
$\cH$ with the opposite multiplication.
The map $d:\cH\to\cH\ot\cH,$  satisfying the conditions (a) - (c)
induces an associative multiplication law on $\cH^\ast$
\begin{eqnarray*}
&{\langle{\alpha\circ\beta,a}\rangle}={\langle{\alpha\ot\beta,d(a)}\rangle}=
{\langle{\alpha \otimes \beta,f_{1}a_{(1)}\otimes f_{2}a_{(2)}}\rangle}=\\[3pt]
&{\langle{\alpha ,f_{1}a_{1}}\rangle}\,{\langle{\beta,f_{2}a_{(2)}}\rangle}=
{\langle{\Delta^{*}(\alpha ),f_{1}\otimes a_{(1)}}\rangle}\,
{\langle{\Delta^{*}(\beta),f_{2}\otimes a_{(2)}}\rangle}=\\[3pt]
&{\langle{\alpha_{(1)},f_{1}}\rangle}\,{\langle{\alpha_{(2)},a_{(1)}}\rangle}\,
{\langle{\beta_{(1)},f_{2}}\rangle}\,{\langle{\beta _{(2)},a_{(2)}}\rangle}= \\[3pt]
&{\langle{\alpha_{(1)},f_{1}}\rangle}\,{\langle{\beta _{(1)},f_{2}}\rangle}\,
{\langle{\alpha_{(2)},a_{(1)}}\rangle}\,{\langle{\beta_{(2)},a_{(2)}}\rangle}=\\[3pt]
&{\langle{\alpha_{(1)},f_{1}}\rangle}\,{\langle{\beta_{(1)},f_{2}}\rangle}\,
{\langle{\alpha_{(2)}\otimes \beta_{(2)},a_{(1)}\otimes a_{(2)}}\rangle}= \\[3pt]
&{\langle{\alpha_{(1)},f_{1}}\rangle}\,{\langle{\beta _{(1)},f_{2}}\rangle}\,
{\langle{\alpha_{(2)}\otimes \beta_{(2)},\Delta(a)}\rangle}=
{\langle{\alpha_{(1)}\otimes\beta_{(1)},\mathcal{F}}\rangle}\,
{\langle{\alpha_{(2)}\otimes\beta_{(2)},\Delta (a)}\rangle}= \\[3pt]
&{\langle{\alpha_{(1)}\otimes\beta_{(1)}\otimes\alpha_{(2)}\otimes\beta_{(2)},
\mathcal{F}\otimes\Delta(a)}\rangle}=
{\langle{\sigma_{23}\left(\alpha_{(1)}\otimes\alpha_{(2)}\otimes\beta_{(1)}
\otimes\beta_{(2)}\right),\mathcal{F}\otimes\Delta(a)}\rangle}=\\[3pt]
&{\langle{\sigma_{23}\Delta^*(\alpha)\otimes\Delta^*(\beta),\mathcal{F}\otimes\Delta(a)}\rangle}\,.
\end{eqnarray*}
Hence $\alpha\circ\beta ={\langle{\alpha_{(1)}\ot\beta_{(1)},{\cF}}\rangle}\alpha_{(2)}\beta_{(2)}$.
\begin{prp}
Twist elements are in a one-to-one correspondence
with module algebra structures on $\cH^\ast$ (over $\cH_{op}$), such that
$\varepsilon \circ f=f\circ\varepsilon=f$\,.
\end{prp}
\begin{cor}
Let $\{e_0=1, e_1, e_2, \ldots \}$  and
$\{e^0=\varepsilon, e^1, e^2, \ldots\}$
be the dual bases in $\cH$ and $\cH^\ast,$ respectively, and let
$$e^i \circ e^j=\sum_k m^{ij}_k e^k.$$
Then
$$
d(1)=\sum_{r,s} m^{rs}_0 e_r \otimes e_s\quad\text{and}\quad
m^{0k}_r=m^{k0}_r=\delta^k_r.
$$
\end{cor}
Indeed,
\begin{eqnarray*}
\left\langle e^j\circ e^i,e_p\right\rangle
&=&\left\langle e^j\otimes e^i,d(e_p)\right\rangle  \\[3pt]
\left\langle e^j\circ e^i, e_p\right\rangle
&=&\left\langle \sum\nolimits_k
\stackrel{\circ }{m}_k^{j\,i}e^k, e_p\right\rangle
=\sum\nolimits_k\stackrel{\circ }{m}_k^{j\,i}\left\langle
e^k, e_p\right\rangle =\stackrel{\circ }{m}_p^{j\,i} \\[3pt]
\left\langle e^j\otimes e^i,d(e_p)\right\rangle  &=&\left\langle e^j\otimes
e^i,\sum\nolimits_{r,s}d_p^{r\,s}e_r\otimes e_s\right\rangle
=\sum\nolimits_{r,s}d_p^{r\,s}\left\langle e^j\otimes e^i,e_r\otimes
e_s\right\rangle = \\[3pt]
&=&\sum\nolimits_{r,s}d_p^{r\,s}\left\langle e^j,e_r\right\rangle
\left\langle e^i,e_s\right\rangle =d_p^{j\,i}.
\end{eqnarray*}
Thus $d_p^{j\,i} ={\stackrel{\circ}{m}}_p^{j\,i} $ and
\begin{eqnarray*}
d(e_p) &:=&\sum\nolimits_{r,s}d_p^{r\,s}e_r\otimes e_s=\sum\nolimits_{r,s}
\stackrel{\circ }{m}_p^{r\,s}e_r\otimes e_s \,,\\[3pt]
d(e_p) &=&d(1e_p)=d(1)\Delta (e_p)=\mathcal{F}\Delta (e_p)\,.
\end{eqnarray*}

\section{Examples}

\subsection{The case of $U(b_2)$}
Let ${\cH}=U(b_2)$ be the universal enveloping algebra of the Borel
algebra with generators $h, x$ and the defining relation  $[h,x]=x$. It is
known, that monomials $e_{m,k}=h^m x^k, m, k=0, 1, 2,\ldots$
form a linear basis in ${\cH}$ (PBW-theorem). Let us consider
the jordanian $R$-matrix
$$
R=\left ( \begin{array}{cccc}
1 & -\xi & \xi & \xi^2\\
0 & 1 & 0 & -\xi\\
0 & 0 & 1 & \xi\\
0 & 0 & 0 & 1
\end{array}\right )
$$
which is a unitary solution of the Yang - Baxter equation~\cite{M}
(that is $R_{21}R=1\otimes 1$) and quasiclassical. Therefore according
to~\cite{D2} there exists a twist $\cF \in {\cH} \otimes {\cH}$  such,
that $R=(\rho \otimes \rho)(\cF_{21}{\cF}^{-1})$,
where $\rho$ is the vector representation of $U(b_2)$.
To find $\cal F$, we shall consider the following algebra,
(it is not a Hopf algebra):
$$
{\cA}=\frac{{\mathbb C}\langle t_{11}, t_{12}\rangle}{R T_1\cdot T_2 =
T_2\cdot T_1}
$$
Here
$$
T_1  = T\otimes 1= \left ( \begin{array}{cccc}
t_{11} & 0 & t_{12} & 0\\
0 & t_{11} & 0 & t_{12}\\
0 & 0 & t^{-1}_{11} & 0\\
0 & 0 & 0 & t^{-1}_{11}
\end{array}\right )
$$
\smallskip
and
\smallskip
$$
T_2  = 1 \otimes T=\left ( \begin{array}{cccc}
t_{11} & t_{12} & 0 & 0\\
0 & t_{11^{-1}} & 0 & 0\\
0 & 0 & t_{11} & t_{12}\\
0 & 0 & 0 & t^{-1}_{11}
\end{array}\right )
$$
\smallskip
\smallskip
Let ${\mathbb C}\langle t_{11}, t_{12} \rangle$ be an algebra,
generated by the noncommutative elements $t_{11}, t_{12}$ satisfying
the relation $RT_1\cdot T_2=T_2\cdot T_1$.
\begin{prp}
\begin{enumerate}
\item[\rm{1)}] $\displaystyle {\cA}=\frac{{\mathbb C}\langle t_{11},
     t_{12}\rangle}{[t_{11}, t_{12}]=\xi \cdot 1}$
\item[\rm{2)}] ${\cA}$ is a left module algebra over
  ${\cH}=U(b_2)$ with the action:
$$
\tilde{x}t_{12}=t_{11},\ \tilde{h}t_{12}=-\frac{1}{2} t_{12},\
\tilde{x}t_{11}=0,\ \tilde{h}t_{11}=\frac{1}{2} t_{11}.
$$
\end{enumerate}
\end{prp}
\begin{proof}
Both statements can be checked immediately. It is also possible to prove
these statements as follows. First of all, it is easy to prove
the Poincare - Birkhoff - Witt theorem: the  monomials
$$
f^{k,m}=t_{11}^k t^n_{12}
$$
form a linear basis in ${\cA}$. Let us consider the algebra
of polynomials on the group $B_2$
$$
Fun = \frac {{\mathbb C} \langle
  t_{11}, t_{12}\rangle}{T_1T_2=T_2T_1}={\mathbb C}[t_{11},t_{12}]\,.
$$
It is known, that $Fun$ is both a left $U(b_2)$ - module algebra and a right
$U(b_2)_{op}$ - module algebra. It is also known that $U(b_2)$ acts on
$Fun$ by  left derivatives whereas $U(b_2)_{op}$ acts by right derivatives.
Because of in the defining relation $R$ stands on the left, it is possible to
show, that the structure of the left $U(b_2)$ - module algebra on ${\cA}$ is
inherited from $Fun$.
\end {proof}
\begin {cor}
As left $U(b_2)$ - module algebras ${\cA}$ and $Fun$ are isomorphic.
\end {cor}
Let us recall, that the used pairing between $U(b_2)$ and $Fun$
is determined by the relation
$\langle f, b\rangle = f\wt{b}(e)$,
in which $f\in Fun, b \in U(b_2),$  and $e$ is the unit element
of the group $B_2$.

So, we have to find dual bases in $U(b_2)$ and $Fun$.
However, it is known, that it is not possible therefore we expand
algebra $Fun$. We will construct such an extension according to the
method suggested in~\cite{KMjgp}.
We recall, that $Fun$ is the Hopf algebra, and the mentioned above
dualization is the pairing of the Hopf algebras $Fun$ and $U(b_2)_{op}$.

Let $\varphi$ and $\omega$ denote the elements dual to $h$ and $x,$
respectively. Then
$$
\langle h^kx^m, \varphi^p\omega^q\rangle =
k!\,m!\cdot \delta_{k,p} \cdot \delta_{m,q},
$$
and the tensor
$$
{\mathcal T}=\exp(h\otimes\varphi)\exp(x\otimes\omega)\in
U(b_2)\otimes\widetilde{Fun}
$$
is the canonical element. Here we have denoted the Hopf algebra, generated
by the elements $\varphi,\omega,e^{\pm{\varphi}/{2}}$ by $\widetilde{Fun}$
with $\Delta \varphi,\, \Delta \omega$ defined by
$$
\langle \Delta(a),h^k x^m \otimes h^p x^q\rangle=
\langle a, h^k x^m h^p x^q\rangle.
$$
It is possible to check, that $\varphi$ is a primitive element, and
$\Delta\omega=\omega\otimes e^{-\varphi}+1\otimes\omega$.
Thus, it is clear that $\widetilde{Fun}$ is a right module algebra
over $U(b_2)_{op}$.

If  ${\mathcal F} \in U(b_2) \otimes U(b_2)$ is a twist, then we will
define a new multiplication in $\widetilde{Fun}$ as follows: since
$\mathcal{F}$ defines a coassociative operation
$d_{\mathcal{F}} (a)=\mathcal{F}\Delta (a)$ on $U(b_2),$ then by virtue of
the duality between $U(b_2)$ and $\widetilde{Fun}$
we obtain the following associative
multiplication on $\widetilde{Fun}$, which can be described as follows.

Let us define
$$
{\mathcal T}_1=\exp (h\otimes 1\otimes\varphi)\exp(x\otimes 1\otimes\omega)
\quad\text{and}\quad
{\mathcal T}_2=1\otimes{\mathcal T},
$$
and also
$$
{\wt{T}}_1= (\rho \otimes \rho\otimes id){\mathcal T}_1
\quad\text{and}\quad
{\wt{T}}_{2}=(\rho \otimes \rho\otimes id){\mathcal T}_2,
$$
where $\rho$ denotes 2-dimensional vector representation of the algebra
$U(b_2)$.
Then  ${\mathcal T}_1{\mathcal T}_2={\mathcal T}_2{\mathcal T}_1.$
Further we will define ${\mathcal T}_1\circ {\mathcal T}_2$  by
$$
{\mathcal T}_1\circ  {\mathcal T}_2=\mathcal{F}_{12}
(\Delta\otimes id)(\mathcal{T}).
$$
Then it is possible to prove, that
$$
{\mathcal T}_2\circ  {\mathcal T}_1=\mathcal{F}_{21}
(\Delta_{op}\otimes id)(\mathcal{T})
$$
Indeed, if we represent $\mathcal{T}$ as $\sum e_i\otimes e^i,$ with
dual bases $\{ e_i \}, \{e^i\}$ in $U(b_2)$  and $\widetilde{Fun},$
respectively, we obtain
$$
{\mathcal T}_1\circ {\mathcal T}_2=
\sum e_i\otimes e_j\otimes e^i\circ e^j=
{\mathcal F}_{12}(\Delta\otimes id)({\mathcal T}).
$$
Calculating ${\mathcal T}_2\circ  {\mathcal T}_1,$  we get
$$
{\mathcal T}_2\circ {\mathcal T}_1=
\sum e_i\otimes e_j\otimes e^j\circ e^i=
{\mathcal F}_{21}(\Delta_{op}\otimes id)({\mathcal T}).
$$
Let ${\mathcal R}={\mathcal F}_{21}{\mathcal F}^{-1}$ and
$R=(\rho \otimes \rho)({\mathcal R})$.  Due to the
equality $\Delta=\Delta_{op},$
we have proved
\begin {thm}
${\mathcal R} {\mathcal T}_1 \circ {\mathcal T}_2
={\mathcal T}_2 \circ {\mathcal T}_1$
and, consequently,
$R\tilde{T}_1\circ\tilde{T}_2=\tilde{T}_2\circ\tilde{T}_1$.
\end {thm}
\vskip0.3cm

This theorem allows us to construct a twist for the given
jordanian $R$-matrix
(and not only for the jordanian one).

Taking into account, that
$$
h=\frac{1}{2}\,\sigma^z=
\frac{1}{2}\,\left (
  \begin{array}{cc}
1 & 0\\
0 & -1
 \end{array}
\right ),\qquad
x=\sigma^+=\left (
\begin{array}{cc}
0 & 1\\
0 & 0
\end{array}
\right ),
$$
after simple calculations  we obtain
$$
\widetilde{T}=\left ( \begin{array}{cc}
e^{\varphi/2} & e^{\varphi/2} \cdot \omega\\
0 & e^{-\varphi/2}
\end{array}\right ).
$$
Then it is possible to rewrite the defining relation
$R\wt{T}_1\circ\wt{T}_2=\wt{T}_2\circ\wt{T}_1$
as follows
$$
e^{\varphi/2} \circ \omega =\omega \circ e^{\varphi/2}+\xi e^{-\varphi/2}\,.
$$
\begin {rem}
Let us recall, that
$$
{\cA}=\frac{{\mathbb C}\langle t_{11}, t_{12}\rangle}{ R T_1\cdot T_2=
T_2\cdot T_1}
$$
is a right module algebra over ${\cH}_{op}=U(b_2)_{op}$. Standard
considerations (note, that we should start from the algebra $\wt{Fun}$,
which is a left ${\cH}$-module and a right ${\cH}_{op}$-module)
show that the algebra, generated by $\wt{T}$ with the defining relation
$R\wt{T}_1\circ\wt{T}_2=\wt{T}_2\circ\wt{T}_1$
is a right ${\cH}_{op}$-module isomorphic to $\wt{Fun}$
as a right ${\cH}_{op}$-module. Advantage of using $\wt{T}$ instead of $T$
is, the existence of the dual bases in this case. The isomorphism
is given by the correspondence:
$$
\varphi^p \omega^q \longrightarrow \varphi \circ \varphi \ldots \circ
\varphi \circ \omega \circ \ldots \circ \omega\,,
$$
where we use a short-hand notation $\varphi^p \circ \omega^q.$
\end {rem}

Let us rewrite the relation
$$
e^{\varphi/2} \circ \omega=\omega \circ e^{\varphi/2}
+ \xi e^{-\varphi/2}
$$
as
$$
\varphi\circ\omega=\omega\circ\varphi +2\xi e^{-\varphi}.
$$
By induction we obtain
$$
e^{n\varphi}\circ\omega -\omega\circ e^{n\varphi}=
2 n\xi e^{(n-1)\varphi}.
$$
Using an analytic continuation (formally it is necessary to introduce
$z=e^\varphi$) we come to conclusion that
$$
e^{a\varphi}\circ \omega - \omega \circ e^{a\varphi}=2a\xi
e^{(a-1)\varphi}\,.
$$

Using induction in $m$ we get
\begin{equation}\label{ast}
\omega^m \circ e^{a\varphi}=
\sum^m_{p=0} (-1)^p \binom{m}{p} a(a-1)\ldots (a-p+1)(2\xi)^p
e^{(a-p)\varphi} \circ \omega^{m-p}\,.
\end{equation}
Recall, that to find the twist it is necessary to compute the
coefficient at 1 in the expression
$(\varphi^n \circ \omega^p) \cdot (\varphi^m \circ \omega^q)$
when we reduce it to the form
$\sum C^{np, mq}_{rs} \varphi^r \circ \omega^s.$
In other words, it is necessary to find $C^{np, mq}_{00}.$

From the relation \reff{ast} it follows that if $n$ and $q$ are
different from zero then $C^{qm,kn}_{00}=0.$
So, it is necessary to compute $C^{0m,k0}_{00}$.
Moreover, it is obvious that it is enough to consider only last
term, namely
$$
(-1)^m (2\xi)^m a(a-1)\ldots (a-m+1)e^{(a-m)\varphi}\,.
$$
To find the commutation relations for the elements
$\omega^m$ and $\varphi^k$ we must differentiate the relation \reff{ast}
in ``variable $a$ " several times and then put $a=0$. In particular, it is
easy to see, that
$$
 C^{0m,k0}_{00}=\frac{(-2\xi)^m}{m!k!} \cdot
\frac{d^k}{da^k} (a)_m \biggl.\biggr|_{a=0}\,,
$$
where $(a)_m=a(a-1)\ldots (a-m+1).$  Then
$$
{\mathcal F}=\sum^\infty_{m=0} \frac{(-1)^m}{m!} (2\xi)^m \cdot x^m
\otimes \sum^\infty_{k=0} \frac{h^k}{k!} \frac{d^k}{da^k}(a)_m
\biggl.\biggr|_{a=0}\,.
$$
where $h,x$ are generators of the algebra $U(b_2)$.

The second factor in the tensor product is the expansion of
$(a)_m$  in the Taylor series in a neighborhood of the point $a=h.$
Thus, finally it yields
$$
{\mathcal F}=\sum^\infty_{m=0} \frac{1}{m!} (-2 \xi x)^m \otimes
(h)_m=(1\otimes 1-2\xi x \otimes 1)^{1\otimes h}
=\exp(\ln(1-2\xi x)\ot h)\,.
$$
It is easy to check that in the fundamental representation the
jordanian $R$-matrix has been reproduced
$$
(\rho \otimes \rho) ({\mathcal F}_{21} {\mathcal F}^{-1})=\left (
\begin{array}{crcr}
1 & -\xi & \xi & \xi^2\\
0 & 1 & 0 & -\xi\\
0 & 0 & 1 & \xi\\
0 & 0 & 0 & 1
\end{array}\right ).
$$

\subsection {One more example of twist for $U(b_2)$}
As another elementary example we will consider the following unitary
quasiclassical solution of the Yang - Baxter equation
(see, for example, \cite{KLMjmp}):
\begin{equation}\label{b1}
{\mbf{R}}_{triang}=
\begin{pmatrix}
1&-\xi&\xi&-{\xi}^2\\
0&1&0&\xi\\
0&0&1&-\xi\\
0&0&0&1
\end{pmatrix}
=1 + \xi r + {\cal O}(\xi^2)\,,
\end{equation}
where $r=\sigma^+\ot\sigma^0-\sigma^0\ot\sigma^+$ and
$\sigma^0=\text{diag}(1,1).$ From the structure of the classical
$r$-matrix it follows that for the
construction of a carrier for the $r$-matrix
it is necessary to add the central
element $c$ to the Borel algebra from the previous example.
As dual bases in $U$ and $U^*$ we will choose
\begin{equation}\label{b2}
h^m\,x^k\,c^p\quad\text{and}\quad
\frac{\varphi^m\,\omega^k\,\alpha^p}{m!\,k!\,p!}
\end{equation}
The canonical element in $U\ot U^*$ is equal to
\begin{equation}\label{b3}
{\cT} = \exp\left(h\ot\varphi\right) \exp\left(x\ot\omega\right)
\exp\left(c\ot\alpha\right),
\end{equation}
and in the representation
$\rho(c)=1,\,\, \rho(h)=\half\sigma^z,\, \rho(x)=\sigma^+$
we have
\begin{equation}\label{b4}
\cT=e^{h\ot\varphi}e^{x\ot\omega}e^{c\ot\alpha}
\biggl.\biggr|_{\{\rho:c=1\,,\, h=\half\sigma^z\,,\, x=\sigma^+\}}=
\begin{pmatrix}e^{\varphi/2}&0 \\ 0&e^{-\varphi/2} \end{pmatrix}
\begin{pmatrix}1&\omega \\ 0&1 \end{pmatrix}
\begin{pmatrix}e^{\alpha}&0 \\ 0&e^{\alpha} \end{pmatrix}
\end{equation}
The relation $RTT=TT$ with the $R$-matrix \reff{b1} defines a new
(noncommutative) multiplication law for the generators
$\varphi,\,\, \omega,\,\, \alpha$ of $U^*$
\begin{equation}\label{CR}
\varphi\alpha=\alpha\varphi,\quad\varphi\omega=\omega\varphi,
\quad\alpha\omega=\omega\alpha+\xi\exp(-\varphi).
\end{equation}
To find the structure constants of the new multiplication law we
take into account that the element $\varphi$ is central and the
commutation relations (\ref{CR})
\begin{equation}\label{}
\vp^m\omega^k\alpha^p\circ\vp^{m_1}\omega^{k_1}\alpha^q=
\varphi^{m+m_1}\omega^{k+k_1}\alpha^{p+q} + \ldots\, .
\end{equation}
To compute the twist $\cF$ we must find
nonzero structure constants $m^{(mkp)(m_1k_1q)}_{(000)}.$
It is possible only if we put
$m = m_1 = 0,\,\, k = q = 0,\,\, p = k_1$.
Using  \reff{CR} we we get for
$(k\leq p):$
\begin{equation}\label{}
\alpha^p \omega^k = \omega^k \alpha^p +
\sum_{n=0}^{k} \omega^n \alpha^{p-k+n} C_k^n \xi^{k-n} C_{p}^{k-n}(k-n)!\,
e^{-(k-n)\varphi}\,.
\end{equation}
Thus, the structure constants (with $k=p, n=0$)
we are looking for are equal to
$$
m^{(00p)(0p0)}_{(000)} =\frac{\xi^p}{p!} ,
$$
and the corresponding universal twist has the following form
\begin{equation}\label{F2}
{\cF}=\sum_p\frac{\xi^p}{p!} c^p \ot x^p
=\exp (\xi c \ot x)\,.
\end{equation}
In the vector representation ${\Bbb C}^2$ $\cF$ is a
$4\times 4$-matrix and the
expression $R=F_{21}F^{-1}$ coincides with the initial $R$-matrix.

\subsection {The case of the Lie superalgebra $sb_2$}
Let us consider the Lie superalgebra $sb_2$,
which also has two
generators: one even $h$ and one odd $x$ with the commutation relations
$$
[h,x] = x , \quad [x,x] = 0 \,.
$$
In the universal enveloping superalgebra $U(sb_2)$ the commutator
$[,]$ is understood as a super-commutator
(${\mathbb Z}_2$-graded commutator):
$$[a,b] = ab -(-1)^{p(a)p(b)} ba\,.$$
We will choose as the unitary quasiclassical solution of the $Z_2$-graded
Yang - Baxter equation~\cite {KS} the following $R$-matrix
\begin{equation}\label{mR3}
R =\left(
\begin{array}{cccc}
1 & \eta & -\eta & 0 \\
0 & 1 & 0 & -\eta \\
0 & 0 & 1 & -\eta  \\
0 & 0 & 0 & 1
\end{array} \right),
\end{equation}
where $\eta$ is an odd Grassmannian parameter: $\eta^2=0$.

The generators of the dual algebra
$U^*(sb_2):\{u, \omega:\langle u, h \rangle=1,\,\,
\langle \omega, x \rangle=1\},\,\, \omega^2 = 0$
are mutually commute. The canonical element
${\cal T} \in U \ot U^*$ is equal
\begin{equation}
{\cal T}=e^{h \ot u}\,e^{x\ot\omega}\,.
\end {equation}

Similar to the first example the superalgebra $sb_2$ has the two-dimensional
representation:
$$
T=\left( \begin{array}{cc}
e^{u/2} & -e^{u/2} \omega \\
0 & e^{-u/2}
\end{array} \right),
$$

We fix the following dual bases in $U$ and $U^*$
$$
h^m x^k \quad\text{and}\quad u^m \omega^k  /m! \,,
$$
where $(m = 0, 1, \cdots ; k = 0, 1).$
The relation $RTT=TT$ with $R$-matrix (\ref{mR3}) defines a
new (noncommutative) multiplication law for the generators $u, \omega$
of the dual superalgebra $U^*$
\begin{equation}\label{CR3}
\omega u = u \omega - 2 \eta e^{-u}, \quad \omega^2 = 0.
\end{equation}
To define the structure constants of this new multiplication law
we take into account commutation relations (\ref{CR3}).
In order to calculate $\cal F$ we need to find the structure constants
$m^{(mk)(m_1k_1)}_{(00)}$ different from zero. For this it is necessary
to set $m = 0,\, k_1= 0,\, m_1 = k = 1.$
As a consequence of (\ref{CR3}), we obtain $m^{(01)(10)}_{(00)}=-2\eta$,
and the corresponding universal twist takes the form
\begin{equation}\label{F3}
{\cal F} = 1 - 2 \eta x \otimes h.
\end{equation}
The universal $R$-matrix
\begin{equation}\label{R3}
{\cal R} = 1 - 2 \eta (x \otimes h - h \otimes x) =
\exp (- 2 \eta (x \otimes h - h \otimes x))
\end{equation}
is linear in the generator $x$ due to nilpotency of $x$ and $\eta$.

\subsection {Twisting of quantum superalgebra $U_q(gl(1|1))$}
The quantum superalgebra $U_q(gl(1|1))$ (see Appendix) is an example
of the case with more complicated $R$-matrix. In the
fundamental two-dimensional representation of the quantum superalgebra
$U_q(gl(1|1))$ this $R$-matrix is equal to
\begin{equation}\label{mR4}
R^{(t)} =\left( \begin{array}{cccc}
q & 0 & 0 & \xi \\
0 & q & \omega & 0 \\
0 & 0 & q^{-1} & 0  \\
0 & 0 & 0 & q^{-1}
\end{array} \right),
\end{equation}
The relation $R^{(t)}TT=TTR$ with the $R^{(t)}$-matrix (\ref{mR4})
defines a new (noncommutative) multiplication law for the elements
(generators) of $U_q(gl(1|1))^*: \alpha, \varphi, b, g$
entering the $T$-matrix. To write down explicit forms
of matrices $T_1$ and $T_2$ it is necessary to include signs
related to the $Z_2$-grading (see Appendix):
\begin{equation}\label{T4}
T_1 = T \ot I = \left( \begin{array}{cccc}
a & 0 & u & 0 \\
0 & a & 0 & -u \\
w & 0 & d & 0  \\
0 & -w & 0 & d
\end{array}\right)\,,\quad
T_2 = I \ot T =\left( \begin{array}{cccc}
a & u & 0 & 0 \\
w & d & 0 & 0 \\
0 & 0 & a & u  \\
0 & 0 & w & d
\end{array}\right).
\end{equation}

In the right hand side of the $RTT=TT$-relation we use the $R$-matrix with
$\xi=0$. The only commutation relation, giving nonzero structure
constant $m^{(j)(k)}_{(0)}$ with nonzero multi-indices needed to calculate
the twist $\cal F$ is
$$
u^2 = \frac{\xi}{q + q^{-1}} d^2 ,
\quad \text{or} \quad
b^2 = \frac{\xi}{q + q^{-1}} 1 + \cdots \,,
$$
where the terms containing generators different from unit are
omitted. Thus,
$m^{(0001)(0001)}_{(0000)}= -\xi /(q + q^{-1})$
and the twist element is
$$
{\cal F} = \exp (-\frac{\xi}{q + q^{-1}} e \ot e)
= 1-\frac{\xi}{q + q^{-1}} e \ot e \,.
$$

\section{Discussion}
In a similar way it is possible to consider the quantum group $U_q(g\ell_n).$
It is known, that the dual Hopf algebra $Fun_q (GL_n)$ which is the algebra
of functions on the quantum group, can be obtained from the relation
$R_q T_1 T_2=T_2T_1 R_q$, where $R_q$ is the Drinfel'd - Jimbo $R$-matrix .
Then the algebra, defined by  $R'_q T_1 T_2=T_2T_1 R_q$, will
be a right module algebra
over $U_q(g\ell_n)_{op}$. Here  $R'_q$ is another $R$-matrix such that
 $R'_q{\cP}$  and  ${\cP}R_q$  have the same spectrum
($\cP$ is the permutation operator in ${\mathbb C}^n\otimes{\mathbb C}^n$).
Then there exists a universal twist
${\mathcal F} \in U_q(g\ell_n)^{\otimes 2},$ such that
$R'_q=(\rho \otimes \rho)({\mathcal F}_{21} {\mathcal R}_q
{\mathcal F}^{-1}),$
where ${\mathcal R}_q$ is the universal $R$-matrix and
$\rho: U_q(g\ell_n)\rightarrow End({\mathbb C}^n)$ is fundamental
vector representation of the Hopf algebra $U_q(g\ell_n)$.

\section*{Appendix}
\subsection*{A. Universal $R$-matrix for $gl_q(1|1)$}

The superalgebra $U_q(gl(1|1))$ is generated by two even $h$, $c$, and two
odd $e$, $f$ elements where $c$ is central, and the remaining generators
satisfy the relations
$$[h,e]=2 e,\quad [h,f]=-2 f,
\quad e^2=f^2=0,\quad [e,f]_+= (q^c-q^{-c})/(q - q^{-1}) $$
The $Z_2$-graded structure of Hopf algebra is defined by
the following coproduct
\begin{align*}
\Delta(h) &= h\ot 1 + 1 \ot h, \\ 
\Delta(c) &= c\ot 1 + 1 \ot c, \\ 
\Delta(e) &= e\ot q^c + 1 \ot e, \\
\Delta(f) &= f\ot 1 + q^{-c} \ot f, 
\end{align*}
and counit $\vep(h)=\vep(c)=\vep(e)=\vep(f)=0$.
The algebra $U_q(gl(1|1))$ is quasitriangular with the universal $R$-matrix
$$
{\cal R} = q^{\frac{1}{2}(c\ot h + h\ot c)}(1- \omega \, e\ot f)
\in U_q(gl(1|1))^{\ot 2},
$$
where $\omega =(q-q^{-1})$. This universal $R$-matrix intertwines the
coproduct with the opposite one.

Let us remark that the superalgebra $gl_q(1|1)$
admits an abelian twist because
it has a two-dimensional classical abelian subalgebra spanned by the
elements $c$ and $h$ with primitive coproduct. There is a special choice
of such a twist element that one of the odd generators,
say $e$, becomes a primitive element. Indeed
$$
\Delta_t(e) = e^{\alpha h \ot c} \Delta(e) e^{-\alpha h \ot c}=
e\ot q^c e^{2\alpha c} + 1\ot e = e\ot 1+ 1\ot e
$$
if we set $\alpha =-{1\over 2}\log(q)$.

Hence, there exists another twist $\cT_1=\exp(\xi e\ot e).$ As a result
of twisting by this element we obtain a new universal $R$-matrix
$$
\wt {\cal R} =
\exp(-\xi e\ot e)e^{\alpha c \ot h}{\cal R} e^{-\alpha h \ot c}
\exp(-\xi e\ot e)\,,
$$
having in the matrix representation an off-diagonal parameter.

\subsection*{B. Dual Hopf algebra $U_q(gl(1|1))^*$}

Generators of the coquasitriangular Hopf algebra dual to $U_q(gl(1|1))$
are defined by nonzero pairings
$$
\langle \vp , h \rangle = 1, \quad
\langle \alpha , c \rangle = 1, \quad
\langle b , e \rangle = 1, \quad
\langle g , f \rangle = 1.
$$
The generators $\vp$ and $\alpha$ are even, whereas $b,g$ are odd.
The canonical element
$$
{\cal T} \in U_q(gl(1|1)) \ot U_q(gl(1|1))^*
$$
is defined by the standard exponentials due to nilpotency of the
generators $e,f$ (and $b,g$)
$$
{\cal T} = \exp (f \ot g) \exp (c \ot \alpha)
\exp (h \ot \vp) \exp (e \ot b),
$$
$$
{\cal T} = (1 + f \ot g) \exp (c \ot \alpha)
\exp (h \ot \vp) (1+e \ot b).
$$
The duality relations allows us to calculate coproduct and to find the
multiplication law for the generators of the dual Hopf algebra
\begin{align*}
\Delta_*(\vp) &= \vp \ot 1 + 1 \ot \vp, \\
\Delta_*(\alpha) &= \alpha \ot 1 + 1 \ot \alpha -
    \frac{2\eta}{\omega} b \ot g, \\
\Delta_*(b) &= b \ot \exp(-2\vp) + 1 \ot b, \\
\Delta_*(g) &= g \ot 1 + \exp(-2\vp) \ot g,
\end{align*}
$$
[\alpha, b]= -\eta \, b,\quad [\alpha, g]= -\eta \ g,
\quad b^2 = g^2=0,\quad [b, g]_+= 0.
$$

To check that the canonical element ${\cal T}$ is a bicharacter
$$
(\Delta \ot 1){\cal T} = {\cal T}_{13} {\cal T}_{23}\,,
\quad (1 \ot \Delta_*){\cal T} = {\cal T}_{12} {\cal T}_{13}\,,
$$
it is necessary (using the coproduct $\Delta$ in the Hopf algebra
$U_q(gl(1|1))$ and dual coproduct $\Delta_*$ in $U_q(gl(1|1))^*$)
to take into account a
particular form of the Baker - Campbell - Hausdorff formula
$$
\exp (t X + Y) = \exp (t X) \exp (\frac{1-e^{-t}}{t} Y), \quad
[X, Y] = Y.
$$
This formula allows us to extract a factor in $(1 \ot \Delta_*){\cal T},$
namely,
\begin{multline*}
\qquad
\exp(c\ot (\alpha\ot 1+1\ot\alpha-\frac{2\eta}{\omega} b\ot g))= \\
=\exp(c\ot\alpha\ot 1)\exp(-\frac{(q^c -q^{-c})}{\omega} b\ot g))
\exp (c \ot 1 \ot \alpha)\,,\qquad
\end{multline*}
due to the commutation relation
$$
[c\ot\alpha\ot 1+c\ot 1\ot\alpha, 1\ot b\ot g] = -2\eta c\ot b\ot g.
$$
In the two-dimensional representation of the Hopf super-algebra
$U_q(gl(1|1))$ the elements of the $2\times 2$-matrix
$T=\left(
\begin{array}{cc}
a     &  u   \\
w     &  d
\end{array}
\right)$
can be expressed in terms of the generators of the dual superalgebra
$U_q(gl(1|1))^*$
$$
a=e^{(\alpha + \vp)}\,,\quad
u=e^{(\alpha + \vp)}b,\,\quad
w=ge^{(\alpha + \vp)}\,,\quad
d=e^{(\alpha - \vp)}+ge^{(\alpha +\vp)}b\,.
$$

\vspace {1.0cm}

{ \flushleft {
Department of Mathematics, \\
University, Sweden \\
E-mail: astolin@math.chalmers.se \\
\bigskip
St. Petersburg Department \\
of V.A.Steklov Mathematical Institute RAS \\
E-mail: kulish@pdmi.ras.ru \\
\bigskip
Defense Engineering Technical University \\
E-mail: evd@pdmi.ras.ru}}

\end {document}